\documentclass{amsart}
\usepackage{mathrsfs}
\usepackage{url}		
\usepackage{comment}
\usepackage{mathtools}
\usepackage{amssymb} 

\newtheorem{theorem}{Theorem}[section]

\theoremstyle{definition}

\theoremstyle{remark}

\numberwithin{equation}{section}

\begin{document}
\title[M Diouf]{Improving the error term in the sieve of  Eratosthenes}


\author{Madieyna Diouf}
\address{Department of Mathematics, Arizona State University,
Tempe, AZ 85282}
\email{mdiouf1@asu.edu}
\thanks{}




\date{09/20/2023}


\begin{abstract}
We have devised an alternative approach to sifting integers in the sieve of Eratosthenes that helps refine the error term. \\
Instead of eliminating all multiples of a prime number $p<z$ in the traditional sieve method, our approach solely eliminates multiples of $p$ that have the minimum prime factor of $p$. \\
By leveraging the density of integers with the least prime factor $p$ in this sieve technique, we obtain a reduced error term and an upper bound of $\pi(x)$ that accurately reflects the prime number theorem.
\end{abstract}
\maketitle
\section{Introduction And Statement of Results}
Sieve theory aims to improve the Sieve of Eratosthenes \cite{Eratosthenes} by addressing limitations related to the M\"obius function and parity problem that hinder accurate estimates. Despite the development of advanced techniques like Brun's sieve \cite{Brun} and Selberg's sieve \cite{Sel}, which has assisted in overcoming some of the obstacles posed by the original sieve, improving the error term in the Sieve of Eratosthenes without straying too far from its core concept remains a challenging task due to our limited knowledge of the complex and random behavior of the M\"obius function. \\
We provide a new approach to understanding the sieve of Eratosthenes and overcoming obstacles created by the M\"obius function and parity problem. Our result enhances the precision of the error term that is currently,
\begin{align}
 E=O\left(2^{\pi(z)}\right).\nonumber
\end{align}
We improve it to:
\begin{align}
 E&=O\left(\sum_{p_i<z}\left\{\frac{x}{p_i}\right\} \prod_{p<p_i}\left(1-\frac{1}{p}\right)\right)=O\left(\pi(z)\right).\nonumber
\end{align}
The best-known upper limit of $\pi(x)$ \cite{Kevin1, Kevin2} using  Eratosthenes's sieve is
\begin{align}
\pi(x)&\ll \frac{x}{\log\log x}.
\end{align}
This outcome in $(1.1)$ is achieved by reducing the value of $z$ to $\log x$. However, it is weaker than the Prime Number Theorem. Nevertheless, it demonstrates the applicability of sieve theory.  To obtain a result that aligns more with the prime number theorem, we give a comparable technique that can  maintain $z$ at $\sqrt x$ and improves the upper bound in $(1.1)$ to 
\begin{align}
\pi(x)&\ll \frac{x}{\log x}.
\end{align}
To arrive at the result expressed in (1.2), we used a sifting method that involved a change from the traditional sieve, combining this with an application of the work of Jared Duker Lichman and Carl Pomerance in the density of the set of integers with least prime factor $p$. 
Consequently, the sieve of Eratosthenes gives 
\begin{align}
S(\mathscr{A, P}, z)&=\sum_{d|P(z)}\mu(d)\left\lfloor\frac{x}{d}\right\rfloor.\nonumber\\
&=x\prod_{p_i<z}\left(1-\frac{1}{p_i}\right)+O\left(2^{\pi(z)}\right).\nonumber
\end{align}
while our result shows that
\begin{align}
S(\mathscr{A, P}, z)&=x\prod_{p_i<z}\left(1-\frac{1}{p_i}\right)+O(\pi(z)).\nonumber
\end{align}
Significant results have been achieved in the field and beyond using applications of sieve theory. Examples include Chen's theorem \cite{Chen}, GYT theorem \cite{Gold}, Green-Tao theorem \cite{Tao} and Zhang-Maynard's result \cite{Zhang, Maynard}, to name a few.
\section{Background}
Consider a finite sequence of integers $\mathscr{A}=\{1, 2, ...,\: x\}$ and a set of primes $\mathscr{P}$. For a given $z\geq 2$, let $P(z)=\prod_{p<z, \: p\in \mathscr{P}} p$. Let $S(\mathscr{A, P}, z)$ denote the number of integers in $\mathscr{A}$ that have no prime divisor $p< z$. \\Say $\chi[statement]=1$ if the statement is true, $0$ if false.
Letting $\mu(d)$ be the M\" obius function, the exact formula of the  Sieve of Eratosthenes by Legendre is 
\begin{align}
\pi(x)-\pi(z)+1=\sum_{d|P(z)}\mu(d)\left\lfloor\frac{x}{d}\right\rfloor. \text{\: \: \: \: \: \: \: } \nonumber
\end{align}
We obtain the Legendre identity
\begin{align}
S(\mathscr{A, P}, z) &=\sum_{d|P(z)}\mu(d)\left\lfloor\frac{x}{d}\right\rfloor. \nonumber\\
S(\mathscr{A, P}, z) &=x\sum_{d|P(z)}\frac{\mu(d)}{d}-\sum_{d|P(z)}\mu(d)\left\{\frac{x}{d}\right\}.\nonumber\\
S(\mathscr{A, P}, z) &=x\sum_{d|P(z)}\frac{\mu(d)}{d}+\sum_{d|P(z)}O(1). \nonumber \\
&=x \prod_{p\leq z}\left(1-\frac{1}{p}\right)+O\left(2^{\omega\left(P(z)\right)}\right). \nonumber\\
S(\mathscr{A, P}, z)&=x \prod_{p\leq z}\left(1-\frac{1}{p}\right)+O\left(2^{\pi(z)}\right).
\end{align}
But $(2.1)$ is of limited use, unless $z$ is very small such as $z=\log x$. We propose the following process that allows the sifting level $z$ to grow with $x$, ideally at an optimal $\sqrt{x}$ rate while considering $\mathscr{A}$ as the entire set of integers in $[1, x]$.
\section{Results}
\subsection{Sieving}
Given the positive integers less than or equal to $x$, \\
1) Remove the multiples of $2$ that have $2$ as the smallest prime factor.\\
2) Remove the multiples of $3$ that have $3$ as the smallest prime factor.\\
3) Remove the multiples of $5$ with $5$ as the smallest prime factor.\\
Continue removing all remaining positive integers less than or equal to $x$ that are multiples of $p_i$ and have $p_i$ as the smallest prime factor, where $p_i$ is the $i$-th prime number, until we reach the largest prime $p_r$ that is less than $z$, where $z=\sqrt{x}$ in the original sieve of Eratosthenes.\\
After completing the process, we obtain a set of positive integers with no multiples of primes less than $z$. These integers are the prime numbers between $z$ and $x$ and the number $1$. So, we obtain the relation
\begin{align}
S(\mathscr{A, P}, z) &=x-\sum_{p_i<z}\sum_{\substack{n\leq x\\
                                                      p_i\mid n\\
                                                       \left(n, \frac{P(z)}{p_ip_{i+1}...p_r}\right)=1}} 1.\nonumber\\
S(\mathscr{A, P}, z) &=x-\mathfrak{S}(x, z) 
\end{align}
where
\begin{align}
\mathfrak{S}(x, z) &=\sum_{p_i<z}\sum_{\substack{n\leq x\\
                                                      p_i\mid n\\
                                                       \left(n, \frac{P(z)}{p_ip_{i+1}...p_r}\right)=1}} 1.\nonumber
\end{align}
It should be noted that the summation
\begin{align}
\mathfrak{I}(x, z)&=\sum_{\substack{n\leq x\\
                             p_i\mid n\\
                              \left(n, \frac{P(z)}{p_ip_{i+1}...p_r}\right)=1}} 1. \nonumber
\end{align}
represents the number of multiples of $p_i$ %
 not exceeding $x$ with least prime factor $p_i$. This value can be estimated with great accuracy. Unlike the traditional sieves, where we remove multiples of primes, and these multiples can overlap, meaning that a multiple of $2$ can also be a multiple of $3$ or a multiple of $p<z$; Here, each multiple of $2$ having the least prime factor as $2$ is distinct from any multiple of $3$ having the least prime factor as $3$ ..., which is different from any multiple of $p<z$ with least prime factor $p$. The following arguments estimate $\mathfrak{I}(x, z)$.\\
Given a prime $p_i<z$, we have $\left \lfloor \frac{x}{p_i} \right\rfloor = \frac{x}{p_i}-\left\{\frac{x}{p_i}\right\}$ or simply $\frac{x}{p_i}+O(1)$. By using the latter estimate, we obtain
\begin{align}
&\frac{x}{p_i}+O(1) \text{\: integers less than or equal to $x$ that are multiples of $p_i$. Among these, }\nonumber\\
&\frac{x}{p_i}\left(1-\frac{1}{p_{i-1}}\right)+O(1) \text{ do not divide $p_{i-1}$. Similarly, } \nonumber\\
&\frac{x}{p_i}\left(1-\frac{1}{p_{i-1}}\right)\left(1-\frac{1}{p_{i-2}}\right)+O(1) \text{ are not divisible by $p_{i-1}$ and not divisible by $p_{i-2}$.} \nonumber\\
&...\nonumber\\
&\frac{x}{p_i}\left(1-\frac{1}{p_{i-1}}\right)\left(1-\frac{1}{p_{i-2}}\right)...\left(1-\frac{1}{p_{1}}\right)+O(1) \text{ are not divisible by any prime $p<p_i$.} \nonumber
\end{align}
Therefore, the number of multiples of $p_i$ less than or equal to $x$, with least prime factor $p_i$  is estimated to:
\begin{align}
\mathfrak{I}(x, z)=\sum_{\substack{n\leq x\\
                             p_i\mid n\\
                              \left(n, \frac{P(z)}{p_ip_{i+1}...p_r}\right)=1}} 1&=\frac{x}{p_i}\prod_{p<p_i}\left(1-\frac{1}{p}\right)+O(1). \\
\sum_{p_i<z}\sum_{\substack{n\leq x\\
                                                      p_i\mid n\\
                                                       \left(n, \frac{P(z)}{p_ip_{i+1}...p_r}\right)=1}} 1
&=\sum_{p_i<z} \frac{x}{p_i}\prod_{p<p_i}\left(1-\frac{1}{p}\right)+\sum_{p_i<z}O(1).\nonumber\\
\mathfrak{S}(x, z) &=\sum_{p_i<z} \frac{x}{p_i}\prod_{p<p_i}\left(1-\frac{1}{p}\right)+O(\pi(z)). 
\end{align}
With $(3.1)$ and $(3.3)$, we have
\begin{align}
S(\mathscr{A, P}, z) &= x-\sum_{p_i<z} \frac{x}{p_i}\prod_{p<p_i}\left(1-\frac{1}{p}\right)+O(\pi(z)).
\end{align}
Jared Duker Lichman and Carl Pomerance have  shown in \cite{Pro} that for primes $p$, $q$ and $r$ (the authors noted that $r>1$ may be any real number), if $g(p)$ is the density of the set of integers with least prime factor $p$, then
\begin{align}
\sum_{p\leq r} g(p)=\sum_{p\leq r}\frac{1}{p} \prod_{q\leq p} \left(1-\frac{1}{q}\right) = 1-\prod_{p\leq r}\left(1-\frac{1}{p}\right).
\end{align}
We  multiply equation $(3.5)$ by $x$ and perform variable changes to obtain
\begin{align}
\sum_{p_i<z}\frac{x}{p_i} \prod_{p<p_i} \left(1-\frac{1}{p}\right) &= x-x\prod_{p_i<z}\left(1-\frac{1}{p_i}\right).
\end{align}
By using $(3.6)$ in equation $(3.4)$, we get
\begin{align}
S(\mathscr{A, P}, z) &= x-\left(x-x\prod_{p_i<z}\left(1-\frac{1}{p_i}\right)\right) + O(\pi(z)).\nonumber \\
S(\mathscr{A, P}, z) &=x\prod_{p_i<z}\left(1-\frac{1}{p_i}\right) + O(\pi(z)).
\end{align}
\section{An alternative proof}
In the following proof, using the same sieve, we will exclusively choose the upper limit $x$ of the set of integers $\mathscr{A}=[1, x]$, as a power of $2$ for technical reasons.\\
Our first method in section $3$ bypasses the use of the M\"obius function, resulting in an error term of $O(\pi(z))$. Let's examine our results when we incorporate the M\"obius function. \\After sieving as illustrated in $(3.1)$, we have
\begin{align}
S(\mathscr{A, P}, z) &= x-\mathfrak{S}(x, z).
\end{align}
Where
\begin{align}
\mathfrak{S}(x, z) 
&=\sum_{p_i<z}\sum_{\substack{n\leq x\\
                                                      p_i\mid n\\
                                                       \left(n, \frac{P(z)}{p_ip_{i+1}...p_r}\right)=1}} 1 \nonumber\\
&=\sum_{p_i<z}\sum_{n\leq x }  \chi [ p_i\mid n ] \chi \left[\left(n, \frac{P(z)}{p_ip_{i+1}...p_r}\right)=1\right].\nonumber\\
&=\sum_{p_i<z}\chi \left[\left(n, \frac{P(z)}{p_ip_{i+1}...p_r}\right)=1\right]\sum_{n\leq x }  \chi [ p_i\mid n ].\nonumber\\
&=\sum_{p_i<z} \sum_{\substack{d\big| \left(n, \frac{P(z)}{p_ip_{i+1}...p_r}\right)}} \mu(d)\sum_{n\leq x }  \chi [ p_i\mid n ].\nonumber\\
&=\sum_{p_i<z} \sum_{d\big| \frac{P(z)}{p_ip_{i+1}...p_r}} \mu(d)\chi [ d\mid n ]\sum_{n\leq x }  \chi [ p_i\mid n ].\nonumber\\
&=\sum_{p_i<z} \sum_{d\big| \frac{P(z)}{p_ip_{i+1}...p_r}}\mu(d)\sum_{n\leq x }  \chi [ p_i\mid n ]\chi [ d\mid n ].\nonumber\\
&=\sum_{p_i<z} \sum_{d\big| \frac{P(z)}{p_ip_{i+1}...p_r}}\mu(d)\sum_{\substack{n\leq x \\  p_i\mid n \\d\mid n }} 1.\nonumber\\
&=\sum_{p_i<z} \sum_{d\big| \frac{P(z)}{p_ip_{i+1}...p_r}}\mu(d)\sum_{\substack{n\leq x \\  [ \: p_i, d\:  ] \mid n }} 1.\nonumber
\end{align}
Where $[x, y]$ is the least comment multiple of the integers $x$ and $y$. 
\begin{align}
\mathfrak{S}(x, z) &=\sum_{p_i<z} \sum_{d\big| \frac{P(z)}{p_ip_{i+1}...p_r}}\mu(d)\left \lfloor \frac{x}{[p_i, d]}\right \rfloor.
\end{align}
Moreover, since 
\begin{align}
d\big | \frac{P(z)}{p_ip_{i+1}...p_r}; \nonumber
\end{align}
Then
\begin{align}
gcd(p_i, d)= (p_i, d) =1.
\end{align}
When using the fact that 
\begin{align}
[p_i, d].(p_i, d)= dp_i,
\end{align}
then $(4.3)$ and $(4.4)$ imply that
\begin{align}
[p_i, d]=dp_i. 
\end{align}
After substituting $(4.5)$ into $(4.2)$, the result is
\begin{align}
\mathfrak{S}(x, z) =\sum_{p_i<z} \sum_{d\big| \frac{P(z)}{p_ip_{i+1}...p_r}}\mu(d)\left \lfloor \frac{x}{dp_i}\right \rfloor.
\end{align}
By inserting $(4.6)$ into $(4.1)$, one obtains
\begin{align}
S(\mathscr{A, P}, z) &= x-\sum_{p_i<z} \sum_{d\big| \frac{P(z)}{p_ip_{i+1}...p_r}}\mu(d)\left \lfloor \frac{x}{dp_i}\right \rfloor.\nonumber\\
&=x-\sum_{p_i<z} \sum_{d\big| \frac{P(z)}{p_ip_{i+1}...p_r}}\mu(d)\left(\frac{x}{dp_i}-\left\{\frac{x}{dp_i}\right\}\right ).\nonumber\\
&=x-\sum_{p_i<z} \sum_{d\big| \frac{P(z)}{p_ip_{i+1}...p_r}}\mu(d)\frac{x}{dp_i}+
\sum_{p_i<z} \sum_{d\big| \frac{P(z)}{p_ip_{i+1}...p_r}}\mu(d)\left\{\frac{x}{dp_i}\right\}.\nonumber\\
&=x-\sum_{p_i<z} \frac{x}{p_i}\sum_{d\big| \frac{P(z)}{p_ip_{i+1}...p_r}}\frac{\mu(d)}{d}+
\sum_{p_i<z} \sum_{d\big| \frac{P(z)}{p_ip_{i+1}...p_r}}\mu(d)\left\{\frac{x}{dp_i}\right\}.\nonumber\\
&=x-\sum_{p_i<z} \frac{x}{p_i}\prod_{p<p_i}\left(1-\frac{1}{p}\right)+
\sum_{p_i<z} \sum_{d\big| \frac{P(z)}{p_ip_{i+1}...p_r}}\mu(d)\left\{\frac{x}{dp_i}\right\}.
\end{align}
We use again the result of J.D Lichman and C.Pomerance in \cite{Pro}, previously mentioned here in $(3.5)$. For primes $p$, $q$ and a real number $r>1$, 
\begin{align}
\sum_{p\leq r}\frac{1}{p} \prod_{q\leq p} \left(1-\frac{1}{q}\right) = 1-\prod_{p\leq r}\left(1-\frac{1}{p}\right).
\end{align}
We  multiply equation $(4.8)$ by $x$ and perform variable changes, to obtain
\begin{align}
\sum_{p_i<z}\frac{x}{p_i} \prod_{p<p_i} \left(1-\frac{1}{p}\right) &= x-x\prod_{p_i<z}\left(1-\frac{1}{p_i}\right).\nonumber\\
x-\sum_{p_i<z} \frac{x}{p_i}\prod_{p<p_i}\left(1-\frac{1}{p}\right)&=x\prod_{p_i<z}\left(1-\frac{1}{p_i}\right).
\end{align}
When placing the value of $(4.9)$ into $(4.7)$ we get
\begin{align}
S(\mathscr{A, P}, z) &=x \prod_{p_i< z}\left(1-\frac{1}{p_i}\right)  +\sum_{p_i<z} \sum_{d\big| \frac{P(z)}{p_ip_{i+1}...p_r}}\mu(d)\left\{\frac{x}{dp_i}\right\}.\nonumber\\
S(\mathscr{A, P}, z)&=x \prod_{p_i< z}\left(1-\frac{1}{p_i}\right)  +O\left(\sum_{p_i<z} \sum_{d\big| \frac{P(z)}{p_ip_{i+1}...p_r}}\mu(d)\left\{\frac{x}{dp_i}\right\}\right).
\end{align}
Observe that since $x=2^k$ where $k$ is an integer greater than $1$, 
then the equation
\begin{align}
O\left(\mu(d)\left\{\frac{x}{dp_i}\right\}\right) &= O\left(\frac{\mu(d)}{d}\left\{\frac{x}{p_i}\right\}\right)\nonumber
\end{align}
holds over the $d$, $p_i$ where $p_i<z$ and $d \big| \frac{P(z)}{p_ip_{i+1}...p_r}$. 
\\Therefore,
\begin{align}
O\left(\sum_{p_i<z} \sum_{d\big| \frac{P(z)}{p_ip_{i+1}...p_r}}\mu(d)\left\{\frac{x}{dp_i}\right\}\right) &= O\left(\sum_{p_i<z} \sum_{d\big| \frac{P(z)}{p_ip_{i+1}...p_r}}\frac{\mu(d)}{d}\left\{\frac{x}{p_i}\right\}\right).\nonumber\\
&=O\left(\sum_{p_i<z}\left\{\frac{x}{p_i}\right\} \sum_{d\big| \frac{P(z)}{p_ip_{i+1}...p_r}}\frac{\mu(d)}{d}\right).\nonumber\\
&=O\left(\sum_{p_i<z}\left\{\frac{x}{p_i}\right\}\prod_{p< p_i}\left(1-\frac{1}{p}\right)\right).
\end{align}
With $(4.10)$ and $(4.11)$, we have
\begin{align}
S(\mathscr{A, P}, z) &=x \prod_{p_i< z}\left(1-\frac{1}{p_i}\right)  +O\left(\sum_{p_i<z}\left\{\frac{x}{p_i}\right\}\prod_{p< p_i}\left(1-\frac{1}{p}\right)\right).\nonumber\\
S(\mathscr{A, P}, z) &=x \prod_{p_i< z}\left(1-\frac{1}{p_i}\right)+O(\pi(z)).
\end{align}
Our methods in sections $3$ and $4$ yield identical outcomes namely $(3.7)$ and $(4.12)$. Notably, these methods vastly improve the error term in the Sieve of Eratosthenes.
\section{Retrieving the upper bound in Chebyshev's inequality}
\text{\: \: \; \: \: \; \; \: \: \: \: \: \: \: \: \: \: \: \: }Set $z=\sqrt{x}$ in this application.
\subsection{Upper bound of $S(\mathscr{A, P}, z)$} \text{\: }\\
We know that 
\begin{align}
\prod_{p< z}\left(1-\frac{1}{p}\right)^{-1}=\prod_{p< z}\sum_{m>0}\frac{1}{p^m}>\sum_{k<z}\frac{1}{k}>\int_{1}^{z} \frac{dx}{x}=\log z.
\end{align}
By $(5.1)$, we obtain
\begin{align}
\prod_{p< z}\left(1-\frac{1}{p}\right) < \frac{1}{\log z}=\frac{1}{\log (\sqrt x)}.
\end{align}
As a result of $(5.2)$ and $(4.12)$, we have
\begin{align}
S(\mathscr{A, P}, z) &<\frac{x}{\log (\sqrt x)}+O\left(\pi(\sqrt{x})\right).\\
S(\mathscr{A, P}, z) &\ll \frac{x}{\log x}.\nonumber
\end{align}
\subsection{Upper bound of $\pi(x)$}
\begin{align}
\pi(x)&=\big(\pi(x)-\pi(z)\big)+\pi(z).\nonumber\\
\pi(x)&\leq S(\mathscr{A, P}, z) +\pi(z).
\end{align}
A combination of $(5.3)$ and $(5.4)$, gives
\begin{align}
\pi(x)&\leq \frac{x}{\log {\sqrt x}}+O\left(\pi(\sqrt{x}) \right).\nonumber\\
\pi(x)&\ll \frac{x}{\log x}.
\end{align}
The sieving technique accurately identifies all relevant multiples and avoids duplication for optimal results.
\text{\: }\\

\end{document}